\documentclass[10pt]{amsart}
\usepackage[cp1251]{inputenc}
\usepackage[english]{babel}
\usepackage{amsmath}
\usepackage{amssymb}
\usepackage{amsfonts}
\usepackage{amscd}
\renewcommand{\hat}{\widehat}

\begin{document}

\title{Slowly Changing Vectors and the Asymptotic Finite-Dimensionality of
an Operator Semigroup}

\author{K. V. Storozhuk}

\address{Sobolev Institute of mathematics SB RAS}
\begin{abstract}
Let $X$ be a Banach space and let $T:X\to X$ be a linear power
bounded operator. Put $X_0=\{ x\in X \ \mid \ T^nx\to 0\}$. We prove
that if $X_0\neq X$ then there exists $\lambda \in
\operatorname{Sp}(T)$ such that, for every $\varepsilon >0$, there
is $x$ such that $\|Tx-\lambda  x\|<\varepsilon $ but
$\|T^nx\|>1-\varepsilon$ for all $n$. The technique we develop
enables us to establish that if $X$ is reflexive and there exists a
compactum $K\subset X$ such that
$\liminf\nolimits_{n\to\infty}\rho\{T^nx, K\}<\alpha (T)<1$ for
every norm-one $x\in X$ then $\operatorname{codim} X_0<\infty$.  The
results hold also for a~ one-parameter semigroup.
\end{abstract}

%\begin{keyword}
% Operator semigroup, asymptotic finite-dimension
%\end{keyword}

\maketitle

\section{Introduction}\markboth{K.V. Storozhuk}{Slowly Changing Vectors}

In this article, $X$ is a complex Banach space, $T:X\to X$ is a
linear operator whose all powers are bounded by a constant
$C<\infty$. We use the notations: $B_X$ is the unit ball in $X$,
$\Gamma $ is the unit circle in $\Bbb C$, $X_0=\{x\in X\mid
T^nx\to_{n\to \infty} 0\}.$

A vector $x$ is called an {\it $\varepsilon $-almost eigenvector}
(or simply an {\it $\varepsilon $-eigenvector}) if there exists
$\lambda\in \Bbb C$ such that $\|Tx-\lambda  x\|<\varepsilon $.
These vectors exist for each $\lambda \in
\operatorname{Sp}(T)\cap\Gamma$. In Section~1, we study $\varepsilon
$-eigenvectors that do not shorten much under the iterations $T^n$.

 {\bf{Definition 1}}. Suppose that $\varepsilon >0$.
Call a vector $x\in X$ $\varepsilon $-{\it slow} if
$$
\exists \lambda \in \Gamma  \mid \|Tx-\lambda  x\|<\varepsilon \quad
\text{and}\quad   \|T^n x\|>1-\varepsilon\ \forall n=0,1,2\dots.
$$
For example, the eigenvectors $x$, $Tx=\lambda  x$, $\lambda
\in\Gamma $, are slow.

{\bf{Example 1}}. $T:l_2\to l_2$ is the right shift,
$T(x_1,x_2,\dots)=(0, x_1,x_2,\dots)$. This is an isometry; hence,
every norm-one $\varepsilon $-eigenvector is $\varepsilon$-slow.

{\bf{Example $1^*$}}. $T:l_2\to l_2$ is the left shift,
$T(x_1,x_2,\dots)=(x_2,x_3,\dots)$. The spectrum includes $\Gamma$
but $T^n x\to 0$ for every $x$ and there are no slow vectors.

If a vector $x$ is $\varepsilon $-slow then $T^n x$ are
$C\varepsilon $-slow for each~$n$, since $T(T^nx)-\lambda  T^n
x=T^n(Tx-\lambda x)$. Thus, the angle between the (complex) lines
$T^nx$ and $T^{n+1}x$ is slow not only for $n=0$ but for all
$n\in\Bbb N$ (i.e., the vector changes slowly under the iterations).

{\bf Remark}. Our terminology is by no means connected with the
terms ``slow vector'' and ``slow variable'' of the classical theory
of dynamical systems. In the title, we call slow vectors slowly
changing.

{\bf Definition 2}. An operator $T$ {\it has slow vectors} if, for
every $\varepsilon >0$, there exist $\varepsilon $-slow vectors. An
operator $T$ {\it has many slow vectors} if $\dim X=\infty$ and, for
every $\varepsilon >0$ and $n<\infty$, there exist $n$-dimensional
subspaces in $X$ whose unit spheres consist of $\varepsilon$-slow
vectors.

If the powers $\|T^n\|$ are bounded  below, i.e., there exists a
number $c$ such that $c \|x\|\leq \|T^nx\|$ for every $x$, then
there are many slow vectors. Indeed, if $\|x\|=1$ and $x$ is a
$c\varepsilon $-eigenvector then $\frac xc$ is $\varepsilon $-slow.

If $X_0=X$ then it is clear that there are no slow vectors. It turns
out that the condition $X_0=X$ is the only obstacle to the existence
of slow vectors; if $\operatorname{codim}X_0=\infty$ then there are
many slow vectors (Theorem 1.1).

In Section~2, slow vectors are used in the study of the asymptotic
properties of $T^n$.

It is known that if there exists an attracting compact set $K$,
i.e., such that
$$
\forall x\in B_X \lim\limits_{n\to\infty}\rho(T^nx, K)=0
$$
then $X=X_0\oplus L$, where $L$ is a finite-dimensional invariant
subspace in $X$. This was proved in [1] for Markov semigroups in
$L_1$. For an arbitrary Banach space, this was established in
[2,\,3]. In [4] it was proved that for the splitting $X=X_0\oplus
L$, $\dim L<\infty$, it suffices that a compact set $K$ attract only
sometimes, i.e.,
$$
\forall x\in B_X\ \liminf\limits_{n\to\infty}\rho(T^nx,K)=0.
$$

A semigroup $T^n:X\to X$ or a one-parameter semigroup $T_t:X\to X$
is called {\it asymptotically finite-dimensional} [5] if
$\operatorname{codim} X_0<\infty$. In [6, 1.3.33] the question is
posed whether a semigroup is asymptotically finite-dimensional if
$$
\forall x\in B_X \liminf\limits_{n\to\infty}\rho(T^nx, K)<\alpha
(T)<1.
                                \eqno(*)
$$
Clearly, the condition $(*)$ follows from the above-listed analogous
conditions.

In Section~2 of the article, we prove that if $X$ is reflexive then
an operator satisfying $(*)$ has few slow vectors and
$\operatorname{codim} X_0<\infty$ (Theorem~2.3).  Thus, we give a
partial positive answer to the question of ~[6]. Theorem~2.3 is
easily generalized to the case of a one-parameter semigroup of
$\{T_t:X\to X, t\geq 0\}$.

Note that it is in the reflexive case that the condition
{$\operatorname{codim} X_0<\infty$} for a bounded semigroup implies
the splitting $X_0\oplus L$ [7].

For nonreflexive $X$, the author does not know the answer to the
question of [6, 1.3.33].

\section{Slow Vectors}

{\bf Theorem 1.1}. \it   Suppose that $X$ is a Banach space, $T:X\to
X$, $\|T^n\|<C$. If $X_0\neq X$ then $T$ has slow vectors. If
$\operatorname{codim} X_0=\infty$ then $T$ has many slow vectors.\rm

{\it Proof}.  Without loss of generality, by passing to the
equivalent norm $\|x\|:=\sup\nolimits_n\{\|T^nx\|\}$, we may assume
that $\|T\|\leq 1.$

The scheme of the proof is as follows: (a) $X_0=0\Rightarrow$ slow
vectors exist; (b) $X_0=0\Rightarrow$ there are many slow vectors;
(c)  $\operatorname{codim} X_0=\infty \Rightarrow$ there are many
slow vectors.

\medskip
(a) Introduce the norm $\|x\|_p:=\lim\limits_{n\to\infty}\|T^nx\|$
on $X$. In this norm, $T$ is an isometry. The norms $\|\ \|_p$ and
$\|\ \|$ need not be equivalent but
$$
\|T^kx\|\sim_{k\to\infty}\|T^kx\|_p
                            \eqno(1.1)
$$
for all $x$. If the powers of $T$ are not bounded below then the
space $(X,\| \|_p)$ is incomplete. Let $\widehat{X}$ be the
completion of $X$ in the norm $\|\ \|_p$. Extend the isometry $T$ of
$(X,\|\ \|_p)$ to $\widehat{X}$ and denote the extension by
$\widehat{T}$. Take $\lambda \in
\operatorname{Sp}(\widehat{T})\cap\Gamma$. Suppose that $\varepsilon
>0$ and $\hat{x}\in\widehat{X}$ is a $\|\ \|_p$-one
$\varepsilon$-eigenvector corresponding to $\lambda $. Involving the
fact that $\widehat{T}:\widehat{X}\to\widehat{X}$ is an isometry, we
get
$$
\|\widehat{T}\hat{x}-\lambda \hat{x}\|_p<\varepsilon \quad
\text{and} \quad \forall n\ \|\widehat{T}^n
\hat{x}\|_p=\|\hat{x}\|_p=1>1-\varepsilon.
$$
The set $X$ is dense in $\widehat{X}$. If a vector $x\in X$ is
sufficiently $\|\ \|_p$-close to $\hat{x}$ then all strict
inequalities of the last formula remain valid. Thus, there exists a
vector $x\in X$ such that
$$
\|T x-\lambda  x\|_p<\varepsilon \quad \text{and}\quad  \forall n\
\|T^n x\|_p>1-\varepsilon.
                                \eqno(1.2)
$$
Thus, $x$ is an $\varepsilon$-slow vector of $T$ in the norm $\|\
\|_p$. Of course, $x$ need not be slow in the initial norm, since
$\|T x-\lambda x\|$ may be large. However, applying the $\|\
\|_p$-isometry $T^k$ to (1.2), we infer
$$
\forall k\  \|T(T^kx)-\lambda {T^k x}\|_p<\varepsilon \quad
\text{and}\quad \forall n\ \|T^n (T^k x)\|_p>1-\varepsilon.
                                \eqno(1.3)
$$

By (1.1), starting from some $k=k_0$, inequalities (1.3) also hold
for the norm $\|\ \|$, i.e., starting from some $k$, the vector
$T^kx$ is as well slow in the initial norm of~$X$.

\medskip
(b) Take a number $\lambda \in\Gamma$ in the spectrum of the
isometry $\widehat{T}$ of $(X,\|\ \|_p)$; $\lambda$ has many
$\varepsilon$-eigenvectors (there exist even infinite-dimensional
spheres of $\varepsilon $-eigenvectors (see [8, Chapter~IV,
Theorems~5.33, 5.9]).

Let $l<\infty$ and let $W$ be an $l$-dimensional subspace in $(X,\|\
\|_p)$ whose $\|\ \|_p$-unit  sphere $S$ consists of $\varepsilon
$-eigenvectors. Perturbing $W$ slightly, we may assume that
$W\subset X$. All vectors in $S$ satisfy (1.2) and (1.3). By (1.1),
for each $x\in S$, all the vectors $T^{k_0}x$ are slow vectors for
the operator $T:X\to X$ starting from some $k_0$. The ellipsoid $S$
is compact; therefore, $k_0$ may chosen common for all $x\in S$.
Thus, $X$ includes $(l-1)$-dimensional ellipsoids of the form
$T^{k_0}(S)$ consisting of small vectors.

\medskip
(c)  Consider the quotient space $X/X_0$. Its elements are
$[x]:=x+X_0$. The norm $\|[x]\|$ is as follows:
$\|[x]\|=\rho(x,X_0)=\inf\{\|x-x_0\|,\ x_0\in X_0\}$. We have
$T(X_0)\subset X_0$; therefore, the operator $[T]:X/X_0\to X/X_0 $
is defined. Clearly, $[T]^n=[T^n]$. It is easy to see that if
$[x]\neq [0]$ then $[T]^n[x]\not\to0$. By (b), $[T]$ has many slow
vectors, i.e., for each $l$, in $X/X_0$, there are $l$-dimensional
ellipsoids of slow vectors $[x]$ for $[T]$ (moreover, we may assume
that these ellipsoids have the form $[S]$, where $S$ is an ellipsoid
in $X$):
 $$
\forall n\geq 0\ \|[T]^n[x]\|>1-\varepsilon \quad \text{and}\quad
\|[T][x]-\lambda  [x]\|<\varepsilon \ \forall [x]\in [S].
$$

Since $T^k x_{0}\to 0$ for all $x_0\in X_0$, we have
$$
\forall n\geq 0\ \|T^n(T^kx)\|>1-\varepsilon \quad \text{and}\quad
\|T(T^kx)-\lambda  T^kx\|<_{k\to\infty}\varepsilon
                                \eqno (1.4)
$$
for every $x\in[x]=x+X_0\in[S]$.

The compactness of $S$ enables us to assert now that, starting from
some $k$, the ellipsoids $T^kS\subset X$ consist of slow vectors.
\qed

\medskip
In what follows, we will need some properties of slow vectors.

Denote by $S_{m,\lambda }$ the operator
$\frac1{m+1}\big(\sum\nolimits_{i=0}^m\frac{T^m}{\lambda ^m}\big)$,
the Cesaro mean of  ${T}/{\lambda}$.

{\bf Lemma 1.2}. \it   Suppose that $\delta>0$, $m\in\Bbb N $. If
$T$ has slow vectors then there exists $\lambda \in \Gamma $ such
that
$$
\exists x\in B_X \mid \|S_{m,\lambda }x-x\|<\delta \quad
\text{and}\quad  \forall n \ \|T^n(S_{m,\lambda }x)\|>1-\delta.
\eqno (1.5)
$$
If $T$ has many slow vectors then there exist subspaces $W\subset X$
of an arbitrarily large dimension whose unit spheres $S$ consist of
vectors $x$ satisfying the inequality of {\rm (1.5)}. \rm

{\it Proof}. Involving a compactness of $\Gamma$, consider
$\lambda\in\Gamma $ to which there correspond slow vectors. If
$\varepsilon :=\varepsilon (\delta,m)$ is very small and $x$ is an
$\varepsilon $-slow vector corresponding to $\lambda $ then
$S_{m,\lambda}(x)\approx x$ and $S_{m,\lambda}(x)$ satisfies (1.5).
The remaining part of the proof is obvious. \qed

{\bf Remark}. Geometrically, Lemma 1.2 means that, for each $m$,
there exist spheres of an arbitrarily large dimension that almost do
not flatten under the mappings $T^n(S_{m,\lambda })$ for any $n$.

\section{ Asymptotic Finite-Dimensionality in the~Reflexive Case}

Throughout the section, we suppose that $T:X\to X$ satisfies $(*)$.
We may assume that $K$ is a~ balanced set.

{\bf Lemma 2.1}. \it  Suppose that $x\in B_X$. For each $k$, there
exist vectors $a_1,\dots, a_k\in K$, numbers $m_1>m_2>\dots>m_k$,
and $t_1,\dots,t_k,$ $|t_i|\leq \alpha ^{i-1}$, such that
$$
\|T^{m_1}x-[t_1T^{m_2}a_1+t_2T^{m_3}a_2+\cdots+t_{k-1}T^{m_k}a_{k-1}
+t_k a_k]\|\leq\alpha ^k.
                                \eqno (2.1)
$$
\rm

{\it Proof}. We  write down some inequalities for $k=1,2,3$. The
first is condition $(*)$, and the validity of each subsequent
inequality is guaranteed by an application $(*)$ to the preceding
inequality multiplied by $\alpha $: $$ \exists n_1\mid
\|T^{n_1}x-t_1a_1\|\leq\alpha, \ |t_1|\leq 1,
$$
$$
\exists n_2\mid \|T^{n_2}(T^{n_1}x-t_1a_1)-t_2a_2\|\leq\alpha ^2,\
|t_2|\leq \alpha,
$$
$$
\exists n_3\mid
\|T^{n_3}(T^{n_2}(T^{n_1}x-t_1a_1)-t_2a_2)-t_3a_3\|\leq\alpha ^3,\
|t_3|\leq \alpha ^2,\dots.
$$
To finish, it remains to remove parenthesis and put
$m_j=n_j+\cdots+n_k$. \qed

The sum of the numbers $|t_i|$ in (2.1) does not exceed
$h:=\sum\nolimits_{i=1}^k \alpha ^i=\frac{1}{1-\alpha }$. Hence, the
convex hull $\widehat{K}$ of $\bigcup\nolimits_{i=0}^\infty T^i(hK)$
attracts $B_X$, i.e.,
$$
\forall x\in B_X\ \forall\varepsilon >0\ \exists n\in\Bbb N\
\exists a\in \widehat{K}\mid \|T^nx- a\|<\varepsilon. \eqno (2.2)
$$

We now show that $T$ cannot act by multiplication by a~ scalar on
the subspaces $X$ whose dimension is rather high.

{\bf{Theorem 2.2}}. \it    $\dim\ker(T-\lambda I)<\infty$ for all
$\lambda \in \Gamma$. \rm

{\it Proof}. Choose a finite $(1-\alpha)$-net of $k$ vectors for $K$
and consider the subspace $Y$ spanned by the net. By the
Kre{\u\i}n--Krasnosel$'$ski{\u\i}--Milman Theorem [9], in every
subspace $Z\subset X$ such that $\dim Z>\dim Y$, there exists
a~norm-one vector~ $z$ such that $\rho(z, Y)=1$. By ~$(*)$,
$\rho(T^nz,Y)<\alpha +(1-\alpha )=1$ for some $n$. Therefore, $Tz$
cannot have the form $\lambda \cdot z$. \qed

{\bf{Theorem 2.3}}. \it Suppose that $X$ is reflexive. They $T$
cannot have many slow vectors and so $\operatorname{codim}
X_0<\infty$. \rm

{\it Proof}.  It suffices to prove that to no $\lambda \in\Gamma$
there correspond many slow vectors. We may assume that $\lambda=1$.

By the Statistical Ergodic Theorem (see, for example, [10, \S\,2]),
the operator means $ S_{m,1}=S_m=\frac{1}{m+1}\biggl(\sum\limits_0^m
T^k\biggr)$ converge to the projection $P$ of $X$ onto $\ker(I-T)$.
By Theorem~2.2,
 $\dim\ker(I-T)<\infty$.

On a compact set $K$, the convergence of $S_m-P$ to zero is uniform
(for example, by Arzel\`a's Theorem). The operators $S_m$ commute
with $T$; therefore, the convergence $(S_m-P)\to 0$ is also uniform
on $\widehat{K}$. Hence, starting from some $m$, $\|S_m(a)-P(a)\|$
is sufficiently small for all $a\in \widehat{K}$, for example, less
than $\frac13$. But by (2.2), for every $x\in B$ and large $n$, the
vector $T^nx$ is close to $\widehat{K}$. Therefore,
 $$
\exists m  \forall x\in B_X\
\|(S_m-P)T^nx\|=\|T^n(S_mx)-Px\|\leq_{n\to\infty}1/3.
$$
This implies, for example, that, under $T^n\circ S_m$, every
$k$-dimensional sphere such that $k>\dim\ker(I-T)$, ``flattens''
three times for large $n$ along some radius $x$ ($x$ must be chosen
so that $Px=0$).

This, by Lemma~1.2 and the remark thereto, means that the number
$\lambda =1$ does not have many slow vectors.

The inequality $\operatorname{codim} X_0<\infty$ follows now from
Theorem~1.1. \qed

\end{document}